\documentclass[10pt,reqno]{amsart}
\usepackage{amsmath,bm}
\usepackage{amsfonts}
\usepackage{amssymb}
\usepackage{amsthm}
\usepackage{amscd, times}
\usepackage{pgfplots}

\def\FunctionF(#1){(-3)*(#1+2)*(#1-1)*(#1-3)}%
\usepackage{amsmath}
\usepackage{amssymb}
\usepackage{amsthm}
\usepackage{graphics}
\usepackage{eucal}
\usepackage{mathrsfs}
\usepackage{tkz-fct}

\usepackage[utf8]{inputenc} 
\usepackage{textcomp} 
\usepackage{graphicx}  
\usepackage{flafter}  
\usepackage{bm}  
\usepackage{tikz}
\usepackage[pdftex,bookmarks,colorlinks,breaklinks]{hyperref}  
\usepackage{color}
\usepackage{memhfixc}  
\usepackage{pdfsync}  
\usepackage{amssymb}
\usepackage{amsfonts}
\theoremstyle{plain}
\newtheorem*{maintheorem*}{Main Theorem}
\newtheorem*{thm*}{Theorem}
\newtheorem*{thma*}{Theorem A}
\newtheorem*{thmaa*}{Theorem A'}
\newtheorem*{thmb*}{Theorem B}
\newtheorem*{thmo*}{Theorem 1.1}
\newtheorem*{thmc*}{Theorem C}
\newtheorem*{thmd*}{Theorem D}
\newtheorem*{thmf*}{Theorem 4.1}
\newtheorem*{remark*}{Remark}
\newtheorem*{conjecture*}{Conjecture}
\newtheorem*{prop*}{Proposition}
\newtheorem*{lem*}{Basic Lemma}
\newtheorem{thm}{Theorem}[section]
\newtheorem{cor}[thm]{Corollary}
\newtheorem{lem}[thm]{Lemma}
\newtheorem{prop}[thm]{Proposition}

\usetikzlibrary{calc, positioning,fit,intersections}

\newtheorem*{proofc*}{Proof of Theorem C}

\newtheorem{conjecture}[thm]{Conjecture}

\begin{document}

\author{ Youssef Lazar}
\email{ylazar77@gmail.com}
\date{}

\title{Simultaneous diophantine approximation for a restricted class of pairs of real numbers}

\maketitle

\begin{abstract} We give sufficient conditions on a pair of real numbers for which  the Littlewood conjecture holds. We use the localization of the roots of a cubic equation with coefficients depending on the diophantine properties for the considered pair $(\alpha,\beta)$. The estimates of the roots rely on the properties of the denominators of the convergents of the continued fraction expansion of $\alpha$ and $\beta$.
  \end{abstract}

\section{Introduction}
\label{intro}

Around 1930, J.E.  Littlewood asked the question whether given any two distinct real numbers $\alpha$ and $ \beta$,  the sequence $(n |\sin( n \alpha) \sin( n \beta)|)_{n \geq 0}$ can take arbitrary small values (\cite{little} Problem 5, p.19). The values of such sequences strongly depend on the diophantine properties of the pair $(\alpha, \beta)$. The problem of Littlewood amounts to solving the following conjecture which is called the Litllewood conjecture in the literature. 
\begin{conjecture}[Littlewood] For any $(\alpha, \beta) \in \mathbb{R}^{2}$, one has
 $$~~~~~~ \liminf_{n\rightarrow \infty} ~~n \| n \alpha \|  \| n \beta \| =0.$$
\end{conjecture}

The conjecture is obviously true if $ \alpha $ or $ \beta $ have unbounded partial quotients in their continued fraction expansion or if $1,\alpha $ and $ \beta $ are linearly dependent over the rationals. Thus, it remains to prove the conjecture when $ \alpha, \beta \in   \mathbb{B}$ where  $\mathbb{B}$ is the set of real numbers which have uniformly bounded partial quotients and when $1,\alpha $ and $ \beta $ are linearly independent over the rationals.  \\

This problem remains unsolved but it has been proved for some specific classes of pairs $(\alpha, \beta)  $. 
\begin{itemize}
\item[•] For the pairs $(\alpha,\beta)$,  where $\alpha$ and $\beta$ belongs to the same cubic field. (Cassels-Swninerton-Dyer, 1955, \cite{cs}). Note that it is believed that cubic irrationals are not badly approximable numbers. 
\item[•] For pairs $(\alpha,\beta)$ of badly approximable numbers, such that $\beta \in B(\alpha)$ where $B(\alpha)$ is a subset of $\mathbb{B}$ with $\dim_{H} B(\alpha) =1$. (Pollington-Velani, 2000, \cite{pv}). This is the first result proving the existence of badly approximable solutions. An effective version of this result has been proved by Bugeaud \cite{b} (2014).

\item[•] The first construction of explicit pairs of badly approximable numbers have been discovered by De Mathan \cite{dm} (2003) and Adamczewski-Bugeaud \cite{ab} (2006). Both results gives examples in terms of their continued fraction expansion. 

\item[•] The set of exceptional pairs which do not satisfy the conjecture has been proved to be of Hausdorff dimension zero (Einsiedler, Katok, Lindenstrauss, 2006, \cite{ekl}). This is the strongest evidence towards the conjecture. The full conjecture is implied by a deep conjecture of Margulis on the distribution of orbits under diagonal flows acting on $\mathrm{SL}_{3}(\mathbb{R})/\mathrm{SL}_{3}(\mathbb{Z}) $. 

\item[•]  Lindenstrauss (\cite{l}) suggested to exploit the result of (Einsiedler, Katok, Lindenstrauss, 2006) in order to find a criterion on a real number $\alpha$ such that any pair $(\alpha,\beta)$ satisfies the conjecture. This leads to the notion of combinatorial entropy of a real number. Lindenstrauss  gives the following criterion, if the entropy of a real number $\alpha$ is positive then any pair $(\alpha,\beta)$ satisfies the Litllewood conjecture. Similarily to the dynamical case, the set of reals with null entropy has Hausdorff dimension zero.
\end{itemize}

\noindent For more precise information about the conjecture one can look at the references \cite{q2}, \cite{s} and \cite{v} for the dynamical point of view.
\noindent \subsection{The main result} We give sufficent conditions which implies the validity of the conjecture for a pair of badly approximable numbers. The conditions are given in terms of the denominators of the convergents of the continued fraction expansion of $\alpha$ and $\beta$.

\begin{thm}\label{main} Let $ (\alpha,\beta)\in \mathbb{B}^{2} $. If one can find a sequence $(\eta_{k})_{k \geq 1}$ of real numbers and a pair of sequences $(n_{k})_{k \geq 1}$, $(m_{k})_{k \geq 1}$ of even positive integers with $0 \leq \eta_{k} < 1/3$ for every  $k$,  such that the two following conditions hold for $ k $ large enough\\
\begin{enumerate}
\item $  q_{m_{k}}(\beta)^{11/12 +\eta_{k}/4} \leq  q_{n_{k}}(\alpha)  \leq q_{m_{k}}(\beta)$.\\

\item $\mathrm{lcm}(q_{n_{k}}(\alpha), q_{m_{k}}(\beta)) \leq  q_{m_{k}}(\beta)^{1+\eta_{k}} $.\\

\end{enumerate}

\noindent Then the pair $ (\alpha, \beta) $ satisfies the Littlewood conjecture.
\end{thm}

\noindent In particular, we get the following immediate consequence of the previous theorem by taking $(\eta_{k})=0$ for all $k$, 

\begin{cor}  Let $ (\alpha,\beta)\in \mathbb{B}^{2} $, if one can find two unbounded sequences $(n_{k})_{k \geq 1}$, $(m_{k})_{k \geq 1}$  and $1 \leqslant \lambda_{k} \leqslant  q_{m_{k}}(\beta)^{1/12}$  sequence of positive integers  such that $ q_{m_{k}}(\beta) =  \lambda_{k} q_{n_{k}}(\alpha)$. Then, 
 the pair $ (\alpha, \beta) $ satisfies the Littlewood conjecture.
\end{cor}

\vspace{0.5cm}
\noindent  \textbf{Remarks.} 
(1) The proof relies on the properties of an auxiliary cubic form $F_{k}(t)$ whose coefficients depend on the diophantine properties of $\alpha$ and $\beta$. The conditions (1) and (2) of the theorem \ref{main} have been set in order to ensure, given an $\varepsilon>0$ ,  the existence of an integer $t_{k}$ which is a  multiple of $\mathrm{lcm}(q_{n_{k}}(\alpha), q_{m_{k}}(\beta)) $ and such that $|F_{k}(t_{k})| \leq \varepsilon$. \\

\noindent (2) We expect that this result can be improved by refining the method for  finding a positive multiple of $\mathrm{lcm}(q_{n_{k}}(\alpha), q_{m_{k}}(\beta))$. Indeed,  we have used the fact the length of the interval $I_{j}$ is strictly larger than $\mathrm{lcm}(q_{n_{k}}(\alpha), q_{m_{k}}(\beta)) $ which is asking too much.  A better estimation of the roots of  the cubic equation $F_{k}(t)=\varepsilon$ may lead to relaxing the constraints.  This can be done by using the algebraic expression of the roots of a cubic instead of the trigonometric one.  \\

\section{Rational lines of approximations and cubic equations}

\subsection{Continued fractions expansion of a real number}

 \noindent Let us recall some of the main properties of the theory of continued fractions we are going to use, see for instance Chapter 3 in \cite{ew} and \cite{s}.

\noindent The error in the approximation of $\beta$ by $c_{n}(\beta)$ is denoted by $e_{n}(\beta)$, i.e. $e_{n}(\beta)= \beta - c_{n}(\beta)$, it is exactly given by 
$$ e_{n}(\beta) =   \dfrac{(-1)^{n}}{\beta_{n+1} q_{n}(\beta)+q_{n+1}(\beta)}$$

\noindent where $\beta_{n+1}$ is such that $\beta =[b_{0}; b_{1}, \ldots, b_{n}, \beta_{n+1}] $ (see Lemma 3 E, \cite{s}).  In particular, $e_{2n}(\beta)$ are all positive and one has  ( \cite{ew} ex. 3.1.5 p. 76 )

\begin{equation}\label{error22}
 \dfrac{1}{2q_{2n+1}(\beta)^{2}} \leq e_{2n}(\beta) \leq   \dfrac{1}{q_{2n}(\beta)^{2}}.
\end{equation}

\subsubsection*{Badly approximable numbers}
The elements of $\mathbb{B}$ are the real numbers with bounded partial quotient in their continued fractions representation. It is not difficult to show that 

$$\mathbb{B} := \{ \beta \in \mathbb{R} ~\vert~ \inf_{q \geqslant 1} q\|q\beta \| > 0 ~\mathrm{for}~\mathrm{every}~q\in \mathbb{N}    \}.$$

\noindent  Applying Dirichlet's theorem to a badly approximable number $\beta$ yields the following estimate   
$$\left|| q \beta\right|| \asymp \dfrac{1}{q}. $$  

\noindent One has for any $\beta \in \mathbb{B}$ that $q_{2n+1}(\beta) \ll q_{2n}(\beta)$ for every $n$ thus using (\ref{error22}) we get the important estimate where the constant involved depends on the maximum value $M$ among the $b_{i}$ for all $i \geq 0$,
\begin{equation}\label{error2m}
  e_{2n}(\beta) \asymp \dfrac{1}{q_{2n}(\beta)^{2}}.
\end{equation}  

\subsection{Small values of the Littlewood cubic form}
Let us fix a pair $ (\alpha, \beta) $ of positive real numbers and let us introduce the ternary cubic form deduced from the conjecture
\begin{equation}
f(x,y,z) := x(\alpha x- y)(\beta x-z).
\end{equation}
This cubic is degenerate in the sense that it is the product of three linear forms. The set of zeros of $f$ is the union of the three planes of respective equations $x=0$, $y=\alpha x$ and $z=\beta x$.  If we assume that $\alpha$
and $\beta$ are two irrational numbers, then no nonzero integral vector could lie in the set $ f=0 $. The Littlewood conjecture amounts to proving that
\begin{equation}
m(f)=\inf_{v\in \mathbb{Z}^{3},v_{1}\neq 0} |f(v)|=0.
\end{equation}
This is equivalent to the assertion that for any real $ \varepsilon>0 $, there exists $(x,y,z) \in \mathbb{Z	}^{3}$,  $ x \neq 0 $ such that
$$  |f(x,y,z)| \leq \varepsilon.$$

\noindent Let us fix an arbitrary small real number $ \varepsilon >0$ and set 
$$\mathcal{D}(\varepsilon):= \left\lbrace (x,y,z) \in \mathbb{R}^{3} : |f(x,y,z)| \leq \varepsilon \right\rbrace. $$
If $ x$ is outside the planes $x=0$, $y=\alpha x$ and $z=\beta x$, then  $\mathcal{D}(\varepsilon)$ is bounded by the two parametric surfaces
$$\alpha x -\dfrac{\varepsilon}{x(\beta x-z)} \leq y \leq  \alpha x +\dfrac{\varepsilon}{x(\beta x-z)}$$
and 
$$\beta x -\dfrac{\varepsilon}{x(\alpha x-y)} \leq z \leq  \beta x +\dfrac{\varepsilon}{x(\alpha x-y)}.$$

\noindent One can check that  $\mathcal{D}(\varepsilon)$ contains the slant asymptote $ \mathbb{R}(1,\alpha,\beta) $ intersection of the two planes $ y=\alpha x $ and $z=\beta x$. The set $\mathcal{D}(\varepsilon)$ is quite complicated since it consists in $ 2^{3}$ connected components; each component is contained in one of the intersection of the half-spaces  bounded by the three planes $ x=0 $, $ y=\alpha x $ and $ z=\beta $. Finding one lattice point in any of these 8 connected component will suffices to ensure the existence of a lattice point in all other components by performing a change of sign on each cooordinate. Usual methods of the geometry of numbers fails to apply in this situation due to the irregular shape of the domain $\mathcal{D}(\varepsilon)$. We are going to look for lattice points lying on rational line segments in $\mathcal{D}(\varepsilon)$.

\subsection{Dirichlet lattice points near $\mathcal{D}(\varepsilon)$} For any arbitrary $\varepsilon > 0$, the line $\mathbb{R}(1,\alpha, \beta)$ lies in the domain $\mathcal{D}(\varepsilon)$, and it may happen that one could find a lattice point on it and solving the problem. But a non-trivial pair must satisfy that $1,\alpha$ and $\beta$ are linearly independent over $\mathbb{Q}$, this breaks the hope of finding a lattice point on $\mathbb{R}(1,\alpha, \beta)$.  Nevertheless, Dirichlet's approximation theorem tells us that it is always possible to find a lattice point lying arbitrarily near the line in question with any level of precision. The crucial issue is that the speed of convergence is not strong enough in order to ensure that such integral vectors are in $\mathcal{D}(\varepsilon)$. More precisely, if we consider a sequence of integers given by $(N_{n})_{1 \leqslant i \leqslant n}$.  The two-dimensional  version of Dirichlet's approximation Theorem tells us that for every $n \geq 1$ there exists an integral vector $M_{n}=(x_{n},y_{n},z_{n}) \in \mathbb{Z}^{3}$ with $1 \leq x_{n} \leq N_{n}$ such that

\begin{equation}\label{dirich1}
\left\lbrace \begin{array}{c}
\displaystyle |\alpha x_{n}-y_{n}| \leqslant N_{n}^{-1/2}\\
 \displaystyle |\beta x_{n}-z_{n}|  \leqslant N_{n}^{-1/2}.
\end{array}\right. 
\end{equation} 

Geometrically, this says that the lattice point $M_{n}=(x_{n},y_{n},z_{n})$  approaches the line $\mathbb{R}(1,\alpha, \beta)$ as $n$ increases with a rate of convergence given by $1/\sqrt{N_{n}}$. 
Applying this to the cubic form $f$ we get that
$$|  f(M_{n}) | = |x_{n}| |\alpha x_{n}-y_{n}| |\beta x_{n}-z_{n}| \leq \dfrac{x_{n}}{N_{n}} \leq 1.$$

\noindent Not surprisingly, Dirichlet's Theorem is not enough for proving the Littlewood conjecture. It can happen that  $ |  f(M_{n}) | \leq \varepsilon $ for some positive integer $ n $, in this case the theorem is proved. Thus from now on, we assume that it is not the case, meaning that for every positive integer $ n $, 
 $$ |  f(M_{n}) | \geq \varepsilon. $$ 
 
\noindent Using (\ref{dirich1}), we obtain 
 $$\varepsilon \leq |  f(M_{n}) | = |x_{n}| |\alpha x_{n}-y_{n}| |\beta x_{n}-z_{n}| \leq \frac{x_{n}}{N_{n}}.$$
 
 \noindent This yields the crucial bound on $x_{n}$,
\begin{equation}\label{xnepsilon}
\varepsilon N_{n} \leq  x_{n} \leq N_{n}.
\end{equation}

\noindent Recall that we are assuming that  $ \alpha $ and  $ \beta $  are badly approximable numbers, which amounts to saying that  $$ C :=\min\{\inf_{q\geq1} q\| q\alpha\|,\inf_{q\geq1}q\| q\beta\|\}>0. $$

\noindent In particular, this information together with (\ref{dirich1}) give the following inequalities

\begin{equation}\label{dirichbad}
\left\lbrace \begin{array}{c}
\displaystyle \frac{C}{x_{n}} < \| \alpha x_{n}\|\leq |\alpha x_{n}-y_{n}| \leq1/\sqrt{N_{n}}\\ \\

 \displaystyle \frac{C}{x_{n}} < \| \beta x_{n}\|\leq |\beta x_{n}-z_{n}| \leq 1/\sqrt{N_{n}}.
\end{array}\right. 
\end{equation}

To sum up, given  an arbitrary positive real number $\varepsilon$ and  positive integers $ N_{1},N_{2},\ldots $ we have found a sequence of lattice points $M_{n}=(x_{n},y_{n},z_{n}) \in \mathbb{Z}^{3}$ with $\varepsilon N_{n} < x_{n} \leq N_{n}$ with $ M_{n} \notin \mathcal{D}(\varepsilon)$ but lying within a distance $1/\sqrt{N_{n}} $ to the line $ \mathbb{R}(1,\alpha, \beta) $.

\subsection{Line approximation} Using the lattice points $ (M_{k})_{k \geqslant 1} $ arising from Dirichlet's theorem, we are going to construct a family of rational  lines which are nearly parallel to the line $ \mathbb{R}(1,\alpha,\beta) $ in $\mathbb{R}^{3}$ and which converges towards it. 

 To do so, we introduce the line passing through $ P_{k} $ and directed by the vector given by $(1,c_{n_{k}}(\alpha),c_{m_{k}}(\beta))  $, that is, 
 $$L^{k}_{\alpha, \beta} = M_{k} + \mathbb{R}(1,c_{n_{k}}(\alpha),c_{m_{k}}(\beta)) $$
 where $ c_{n_{k}}(\alpha) $ (resp. $ c_{n_{k}}(\beta) $) is the $n_{k}^{\mathrm{th}}$ (resp. $m_{k}^{\mathrm{th}}$) convergent of $\alpha$ (resp. $\beta $).

\noindent A parametrization of the line segment given by  $v_{k}: \mathbb{R} \rightarrow \mathbb{R}^{3}$ is defined as
$$ v_{k}(t) := M_{k} - t (1,c_{n_{k}}(\alpha),c_{m_{k}}(\beta))=(x_{k}(t), y_{k}(t), z_{k}(t))~~$$ where 
$$(L^{k}_{\alpha, \beta}) :   \left\lbrace \begin{array}{ccc}
x_{k}(t)&=  & x_{k}-t \\
 y_{k}(t)&=  & y_{k}- t c_{n_{k}}(\alpha)\\
 z_{k}(t)&=  & z_{k}- t c_{m_{k}}(\beta).
\end{array}\right.$$

\noindent The strategy for proving  the main theorem consists in showing the existence of an integer $t_{k}$ satisfying the three conditions

\begin{center}
\begin{enumerate}
\item  (Non vanishing condition)  $x_{k}(t_{k}) \neq 0$.
\item (Geometric condition) For every $\varepsilon > 0$, $v_{k}(t_{k}) \in \mathcal{D}(\varepsilon)$ i.e.  $0 < |f(v_{k}(t_{k}))| \leq \varepsilon$. 
\item (Arithmetic conditon) $t_{n}$ is an integral multiple of $l_{k} :=\mathrm{lcm}(q_{n_{k}}(\alpha), q_{m_{k}}(\beta))$.  
\end{enumerate}
\end{center} 

\noindent We have the following obvious bound on $l_{k}$
\begin{equation}\label{lambda1}
 \max\{q_{n_{k}}(\alpha),q_{m_{k}}(\beta) \} \leq l_{k} \leq q_{n_{k}}(\alpha)q_{m_{k}}(\beta).
\end{equation}

\noindent Since the denominators for the convergents of a badly approximable number have a geometrical growth, we can assume that there exists an integer $n_{0}$ large enough such that $q_{n_{k}}(\alpha) \leq q_{m_{k}}(\beta)$ for every $n \geq n_{0}$. In other words, for $k \geq k_{0}$ (\ref{lambda1}) is just
\begin{equation}\label{lambda}
q_{m_{k}}(\beta) \leq  l_{k}  \leq q_{n_{k}}(\alpha)q_{m_{k}}(\beta).
\end{equation}

\noindent Let us introduce the \textit{logarithmic ratio} between $q_{n_{k}}(\alpha)  $ and $ q_{m_{k}}(\beta) $ as $$\gamma_{k} := \dfrac{\ln q_{n_{k}}(\alpha)}{\ln q_{m_{k}}(\beta)}. $$ 
Thus one has  $0<\gamma_{k} \leqslant 1$ for every $ k $ greater than $ k_{0} $, and we can now write 
\begin{equation}\label{rho}
q_{n_{k}}(\alpha) = q_{m_{k}}(\beta)^{\gamma_{k}}.
\end{equation}

\noindent For $n$ large enough, the bounds (\ref{lambda}) can be rewritten 
\begin{equation}\label{lambdaplus}
q_{m_{k}}(\beta) \leq  l_{k}  \leq q_{m_{k}}(\beta)^{1+\gamma_{k}}.
\end{equation}

\noindent The assumptions of Theorem \ref{main} tell us that we have the bounds
\begin{equation}\label{gammabound}
\dfrac{11}{12}+\dfrac{\eta_{k}}{4} <\gamma_{k} \leq 1
\end{equation}
under the condition  $ 0 \leq \eta_{k} <1/3 $ when $k \geq k_{0}$.

\subsection{The cubic polynomial $ F_{k}(t)$ associated to the Dirichlet  lattice points $(M_{k})_{k \geqslant 1}$}

In order to study the existence of lattice points in $\mathcal{D}(\varepsilon)$, we will restrict to study their eventual presence on the lines $(L_{\alpha,\beta}^{k})=\{v_{k}(t) ~:~ t \in \mathbb{R} \}$. First, we evaluate the cubic $f$ on $(L_{\alpha,\beta}^{k})$.  One has, for any real $t$, 

$$f(v_{k}(t)) = (x_{k} - t )\left\lbrace\alpha(x_{k}-t) -(y_{k}-t c_{n_{k}}(\alpha)\right\rbrace \left\lbrace\beta(x_{k}-t) -(z_{k} -tc_{m_{k}}(\beta))\right\rbrace$$
$$ = (x_{k} - t )\left\lbrace (\alpha x_{k}-y_{k})-t e_{n_{k}}(\alpha)\right\rbrace \left\lbrace(\beta x_{k}-z_{k})-t e_{m_{k}}(\beta)\right\rbrace$$

\noindent Let us denote $\delta_{k}(\alpha) :=\alpha x_{k}-y_{k}$ and $\delta_{k}(\beta):=\beta x_{k}-z_{k}$, without loss of generality and changing signs of $y_{k}$ and $z_{k}$ or of $f$ if necessary,  we can assume that $\delta_{k}(\alpha)$ and $\delta_{k}(\beta)$ are both positive. Thus,  we get

$$f(v_{k}(t)) = (x_{k} - t )\left\lbrace \delta_{k}(\alpha)-t e_{n_{k}}(\alpha)\right\rbrace \left\lbrace\delta_{k}(\beta)-t e_{n_{k}}(\beta)\right\rbrace$$
\noindent Also we set  $t_{k}(\alpha)=\delta_{k}(\alpha) /e_{n_{k}}(\alpha)$, $t_{k}(\beta)=\delta_{k}(\beta) /e_{n_{k}}(\beta)$ and $A_{k}=e_{n_{k}}(\alpha) e_{m_{k}}(\beta)$. Hence,

$$f(v_{k}(t)) = -A_{k} (t-x_{k})(t-t_{k}(\alpha))(t-t_{k}(\beta)).$$

\noindent The expression of $ f(v_{k}(t)) $ defines a cubic polynomial in $t$ which we denote $F_{k}(t)$. It assumes three real (positive)\footnote{The positivity of roots is not important but it allows us to work without taking care of the sign in the estimates and inequalities. At the end we are going to solve $|f(v)| < \varepsilon$ with the asbolute value.} roots, namely $x_{k}$, $t_{k}(\alpha)$ and $t_{k}(\beta)$. \\
\begin{figure}[htbp]
\label{dessin2}
\begin{center}

\begin{tikzpicture}
\draw[->] (-1,0) -- (11,0);
\draw (11,0) node[right] {$x$};
\draw [->] (0,0) -- (0,8);
\draw (0,8) node[above] {$y$};

\draw (0,0) node[below] {$0$};

\draw[domain=2:8,scale=1, blue] plot [smooth](\x,\x -15/\x^2);
\draw[domain=1.3:8,scale=1, blue] plot[smooth] (\x,\x +15/\x^2);

\draw [dashed] (7.2,7.2) -- (7.2,0) node[below] {$N_{k}$};
\draw [dashed] (6.6,5.9) -- (6.6,0) node[below] {$x_{k}$};
\draw [dashed] (4.6,3.9) -- (4.6,0) node[below] {$x(t_{k}(\varepsilon))$};

\draw[blue]   (5,2) node[above] {$\displaystyle y= \alpha x - \dfrac{\varepsilon}{\beta x^2}$};
\draw[blue]  (4.5,7) node[above] {$\displaystyle y= \alpha x + \dfrac{\varepsilon}{\beta x^2}$};

\draw  (9,8) node[above] {$\displaystyle y= \alpha x                                                                                   $};

\draw  (3,2) node[above] {$(L^{n}_{\alpha, \beta})$};


\draw[dotted] (7.4,5.6) node[above] {$\bullet(x_{k},y_{k},0)$};
\draw[dotted] (7,7.2) node[right] {$\bullet (N_{k},\alpha N_{k},0)$};



\draw[domain=0:8,scale=1] plot (\x,{\x});
\draw[domain=0.5:8,scale=1,red] plot (\x,{\x -0.7});


\end{tikzpicture}
\caption{A view of the picture in the plane $ z=0 $. The red line $(L^{k}_{\alpha, \beta})$ passes through the Dirichlet lattice vector and cuts  $\mathcal{D}(\varepsilon)$ at $x(t_{k}(\varepsilon))$. The boundary of $\mathcal{D}(\varepsilon)$ is colored in blue. }
\label{default}
\end{center}
\end{figure}
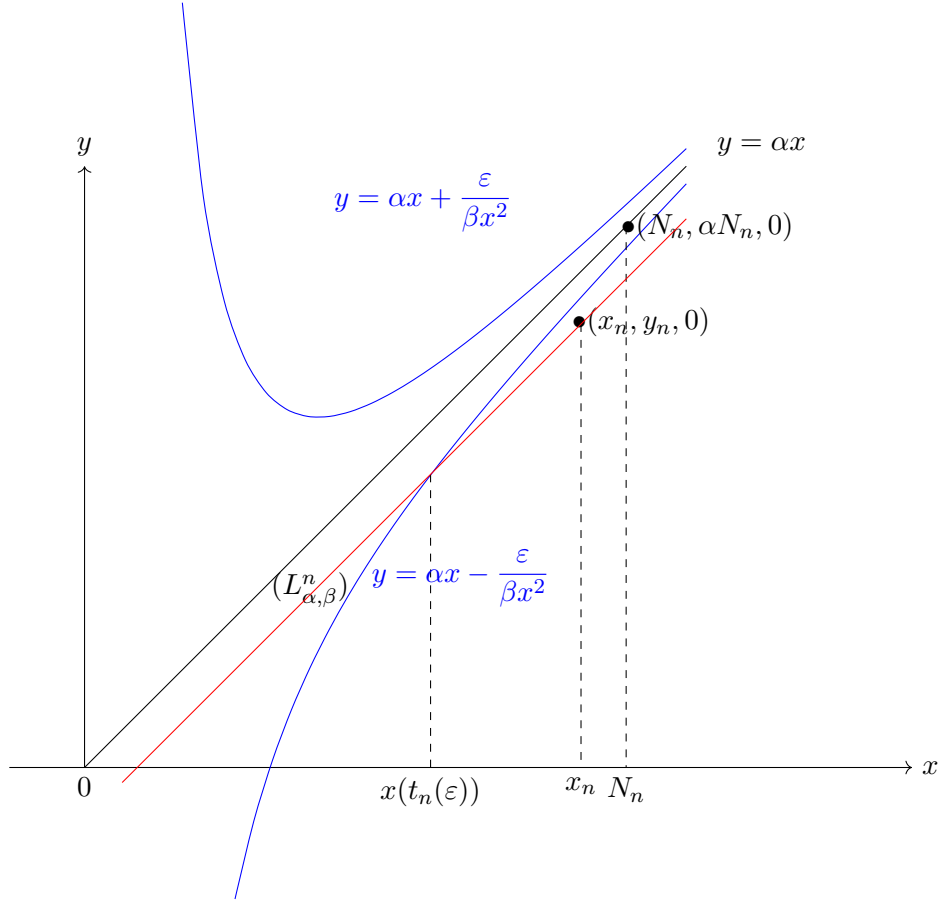

\subsection{Estimates of the coefficients of $ F_{k} $} Let us introduce the elementary symmetric polynomials of the roots $x_{k}$, $t_{k}(\alpha)$ and $t_{n}(\beta)$ of $F_{k}$.
\begin{equation}\label{s}
\left\lbrace \begin{array}{ccc}
  \Sigma_{1}  & = &  x_{k} + t_{k}(\alpha) + t_{k}(\beta)\\
 \Sigma_{2}   & = &  x_{k} t_{k}(\alpha)+  x_{k}t_{k}(\beta) +t_{k}(\alpha) t_{k}(\beta)\\
 \Sigma_{3}   & = &  x_{k} t_{k}(\alpha) t_{k}(\beta).
\end{array}\right. 
\end{equation}

\noindent For the two sequences of real numbers $(\delta_{k})_{k \geqslant 1}$  and $(\gamma_{k})$ define
\begin{equation}\label{Ndelta}
N_{k} =  q_{m_{k}}(\beta)^{\delta_{k}}.
\end{equation}

 \noindent  The following proposition gives asymptotic bounds depending on $\varepsilon$ for each $\Sigma_{i}$ ($i=1,2,3$).
\begin{prop}\label{sigmas} For $n$ large enough, one has
\begin{enumerate}
\item  $  q_{2n}(\beta)^{\delta_{n}} + q_{2n}(\beta)^{2-\delta_{n}} \ll_{\varepsilon} \Sigma_1  \ll_{\varepsilon} q_{2n}(\beta)^{\delta_{n}} + q_{2n}(\beta)^{2-\delta_{n}/2}.$  \\

\item  $q_{2n}(\beta)^{2}  + q_{2n}(\beta)^{2(1+\gamma_{n})-\delta_{n}} \ll_{\varepsilon} \Sigma_{2} \ll_{\varepsilon}  q_{2n}(\beta)^{2(1+\gamma_{n})-\delta_{n}} + q_{2n}(\beta)^{2+\delta_{n}/2}$.\\

\item $ \Sigma_{3} \asymp_{\varepsilon}  q_{2n}(\beta)^{2(1+\gamma_{n})}. $  \\

\end{enumerate}

\end{prop}
\noindent \textit{Proof of the proposition.} (1) From \ref{error22}, \ref{xnepsilon} and \ref{dirichbad} we get
\begin{equation}\label{xn1} 
\varepsilon N_{k} < x_{k} \leq N_{k}, 
\end{equation} 

\begin{equation}\label{tnbeta}
   \dfrac{q_{m_{k}}(\beta)^2}{x_{k}}  \ll    t_k(\beta) \ll \dfrac{ q_{m_{k}}(\beta)^2}{\sqrt{N_{k}}},
\end{equation}

\begin{equation}\label{tnalpha}
  \dfrac{q_{n_{k}}(\alpha)^{2}}{x_{k}}  \ll    t_k(\alpha) \ll \dfrac{q_{n_{k}}(\alpha)^{2}}{\sqrt{N_{k}}}.
\end{equation}

\noindent An estimate of $ \Sigma_{1}$ is thus given by,

\begin{equation}\label{sigma1}
    N_k + \dfrac{1}{N_{k}}\left(q_{n_{k}}(\alpha)^{2}+ q_{m_{k}}(\beta)^{2}\right) \ll_{\varepsilon} \Sigma_{1} \ll_{\varepsilon} N_k + \dfrac{1}{\sqrt{N_{k}}}\left(q_{n_{k}}(\alpha)^{2}+ q_{m_{k}}(\beta)^{2}\right).
\end{equation}

\noindent Using the notation \ref{Ndelta}, 

\begin{equation}\label{sigma11}
q_{m_{k}}(\beta)^{\delta_{k}} + q_{n_{k}}(\beta)^{2\gamma_{k}-\delta_{k}}   + q_{m_{k}}(\beta)^{2-\delta_{k}}  \ll_{\varepsilon} \Sigma_1  \ll_{\varepsilon} q_{m_{k}}(\beta)^{\delta_{k}} + q_{m_{k}}(\beta)^{2\gamma_{k} -\delta_{k}/2}.+ q_{m_{k}}(\beta)^{2-\delta_{k}/2}. 
\end{equation}

\noindent  Since $0 < \gamma_{k} \leq 1$ then, we are reduced to
\begin{equation}\label{sigma12}
  q_{m_{k}}(\beta)^{\delta_{k}} + q_{m_{k}}(\beta)^{2-\delta_{k}} \ll_{\varepsilon} \Sigma_1  \ll_{\varepsilon} q_{m_{k}}(\beta)^{\delta_{k}} + q_{m_{k}}(\beta)^{2-\delta_{k}/2}.
\end{equation}

\noindent (2) Let us consider  $\Sigma_2$. The asymptotic estimates (\ref{tnbeta}), (\ref{tnalpha}) allow us to obtain
\begin{equation}\label{xntnalpha}
q_{n_{k}}(\alpha)^2  \ll x_{k} t_{k}(\alpha) \ll \sqrt{N_{k}} ~q_{n_{k}}(\alpha)^{2},
\end{equation}

\begin{equation}\label{xntnbeta}
q_{m_{k}}(\beta)^2 \ll x_{k} t_{k}(\beta) \ll \sqrt{N_{k}} ~q_{m_{k}}(\beta)^{2}.
\end{equation}

\noindent For the product $t_k(\alpha) t_k(\beta) $, the estimates (\ref{tnbeta}) and (\ref{tnalpha}) do not provide optimal bounds. Instead we use the  assumption that $ |f(x_{k},y_{k},z_{k})| > \varepsilon$, which we assume to hold, otherwise, the conjecture is proved. This condition is equivalent to 
\begin{equation}
\varepsilon < x_{k} \delta_{k}(\alpha) \delta_{k}(\beta)
\end{equation}
\noindent or also,
\begin{equation}
 \dfrac{\varepsilon}{A_{k}} < x_{k}t_{k}(\alpha) t_{k}(\beta).
\end{equation}

\noindent Using (\ref{error2m}) again, we arrive to

\begin{equation}
 \dfrac{\varepsilon q_{n_{k}}(\alpha)^{2}  q_{m_{k}}(\beta)^{2}  }{x_{k}} \ll  t_{k}(\alpha) t_{k}(\beta).
\end{equation}

\noindent Thus, we obtain the following refined estimate

\begin{equation}\label{tnalphatnbeta}
  \dfrac{ q_{n_{k}}(\alpha)^{2}  q_{m_{k}}(\beta)^{2}    }{x_{k}} \ll_{\varepsilon}  t_k(\alpha) t_k(\beta) \ll_{\varepsilon} \dfrac{q_{n_{k}}(\alpha)^{2} q_{m_{k}}(\beta)^{2}}{N_{k}}.
\end{equation}

\noindent From (\ref{xntnalpha}) (\ref{xntnbeta}), (\ref{tnalphatnbeta}),

\begin{equation}\label{sigma21}
   \dfrac{  q_{n_{k}}(\alpha)^{2} q_{m_{k}}(\beta)^{2}}{N_{k}} +\left(q_{n_{k}}(\alpha)^{2}+ q_{m_{k}}(\beta)^{2}\right) \ll_{\varepsilon} \Sigma_{2} \ll_{\varepsilon}   \dfrac{1}{N_{k}} q_{n_{k}}(\alpha)^{2} q_{m_{k}}(\beta)^{2}+  \sqrt{N_k}\left(q_{n_{k}}(\alpha)^{2}+ q_{m_{k}}(\beta)^{2}\right).
\end{equation}

\noindent As $N_{n}=q_{m_{k}}(\beta)^{\delta_{k}}$, we arrive to

\begin{equation}\label{sigma22}
  q_{m_{k}}(\beta)^{2(1+\gamma_{k})-\delta_{k}}  + q_{m_{k}}(\beta)^{2\gamma_{k}} +  q_{m_{k}}(\beta)^{2} \ll_{\varepsilon} \Sigma_{2} \ll_{\varepsilon} q_{m_{k}}(\beta)^{2(1+\gamma_{k})-\delta_{k}}   + q_{m_{k}}(\beta)^{2\gamma_{k} +\delta_{k}/2} + q_{m_{k}}(\beta)^{2+\delta_{k}/2}.
\end{equation}
Thus,

\begin{equation}\label{sigma222}
q_{m_{k}}(\beta)^{2}  + q_{m_{k}}(\beta)^{2(1+\gamma_{k})-\delta_{k}} \ll_{\varepsilon} \Sigma_{2} \ll_{\varepsilon}  q_{m_{k}}(\beta)^{2(1+\gamma_{k})-\delta_{k}} + q_{m_{k}}(\beta)^{2+\delta_{k}/2}.
\end{equation}
and this shows the second assertion of Proposition \ref{sigmas}.\\

\noindent (3) For $\Sigma_{3}$, we already have from (\ref{tnalphatnbeta}) that
\begin{equation}
 \dfrac{\varepsilon}{A_{k}} < x_{k}t_{k}(\alpha) t_{k}(\beta)  \leq q_{n_{k}}(\alpha)^{2} q_{m_{k}}(\beta)^{2}.
\end{equation}

\noindent That is, 

\begin{equation}\label{sigma31}
x_k  t_k(\alpha) t_k(\beta) \asymp_{\varepsilon} q_{n_{k}}(\alpha)^{2} q_{m_{k}}(\beta)^{2}.
\end{equation}

\noindent Hence, 

\begin{equation}\label{sigma3}
 \Sigma_{3} \asymp_{\varepsilon}  q_{m_{k}}(\beta)^{2(1+\gamma_{k})}.
\end{equation}

\noindent  

\begin{flushright}
$  \square$
\end{flushright}

\subsection{Asymptotic shape of the graph $y=F_{k}(t)$} The derivative of $F_{k}(t)$ with respect to the real variable is
$$F^{\prime}_{k}(t)  = -A_{k} \left\lbrace 3t^{2}- 2 \Sigma_{1} t + \Sigma_{2}  \right\rbrace.$$

\noindent To find the local extrema of $F_{n}$ we are then reduced to compute the zeroes of $ 3t^{2}- 2 \Sigma_{1} t +  \Sigma_{2}  $. The fact that we have three roots for $ F_{n} $, a graphical argument shows one has two extrema so  that the discriminant of $ F_{k}^{\prime} $ is nonnegative for every integer $ n\geq 1 $,
 $$\Delta_{k}= 4 \Sigma_{1} ^{2}-12 \Sigma_{2} = 4( \Sigma_{1} ^{2}-3 \Sigma_{2} ).$$

\noindent The roots of the derivative $F^{\prime}_{k}(t)$ are then given by 
\begin{equation}
\tau_{k}^{\pm}=\Sigma_{1} \pm \dfrac{1}{2} \sqrt{\Delta_{k}}.
\end{equation}

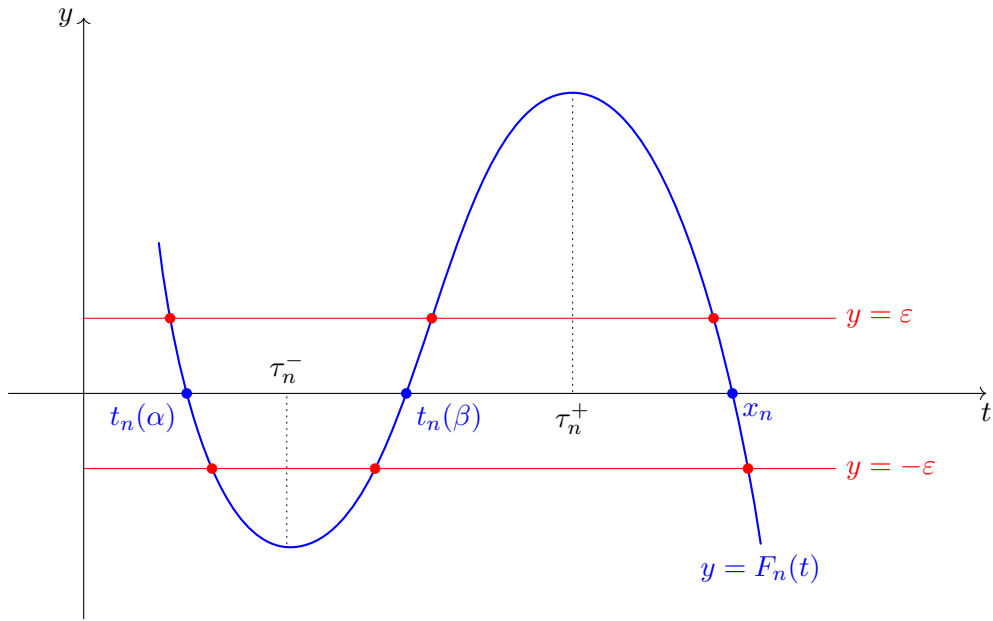
\begin{figure}

\begin{tikzpicture}[scale=1,cap=round]
 \tikzset{axes/.style={}}
 \begin{scope}[style=axes]
 \draw[->,name path=c4] (-1,0) -- (12,0) node[below] {$t$};
 \draw[->] (0,-3)-- (0,5) node[left] {$y$};
 \draw[red,name path=c1] (0,1)-- (10,1) node[right] {$y=\varepsilon$};
 \draw[red,-,name path=c2] (0,-1)-- (10,-1) node[right] {$y=-\varepsilon$};
 \draw [dotted]  (2.7,-2) -- (2.7,0) node[above] {$\tau_{k}^{-}$};
 \draw [dotted]  (6.5,4) -- (6.5,0) node[below] {$\tau_{k}^{+}$};

 \draw [blue,thick,-,name path=c3] plot [smooth,tension=1] coordinates { 
   (1,2) (3,-2) (6.5,4) (9,-2)}  node[below] {$y=F_{k}(t)$}  ;
 
 \fill[red,name intersections={of=c1 and c3}]
    (intersection-1) circle (2pt)
    (intersection-2) circle (2pt)
        (intersection-3) circle (2pt) ;

        \fill[red,name intersections={of=c2 and c3}]
    (intersection-1) circle (2pt)
    (intersection-2) circle (2pt)
        (intersection-3) circle (2pt) ;

          \fill[blue,name intersections={of=c4 and c3}]
    (intersection-1) circle (2pt)node[below left]{$t_k(\alpha)$}
    (intersection-3) circle (2pt)node[below right]{$x_{k}$}
        (intersection-2) circle (2pt) node[below right]{$t_k(\beta)$};


\end{scope}
\end{tikzpicture}
\caption{The shape of $y=F_{k}(t)$}
\end{figure}

\noindent The values of the local maximum and minimum are respectively given by

$$F_{k}(\tau_{n}^{\pm}) = -A_{k} (\tau_{k}^{\pm}-x_{k})(\tau_{k}^{\pm}-t_{k}(\alpha))(\tau_{k}^{\pm}-t_{k}(\beta)).$$

 In fact we have much more precise than this, it requires the use of a result due to Laguerre-Ces\'{a}ro (see e.g. Theorem 6.5.1 \cite{rs}) which says in our case that $\tau^{\pm}_{k}$ is in the middle third of the interval bounded by two the critical numbers bounding  $\tau^{\pm}_{k}$. In concrete words, if the root are ordered as $\tau_{k}^{(1)}<\tau_{k}^{(2)}<\tau_{k}^{(3)} $ so that we have the following configuration
$$ \tau_{k}^{(1)}< \tau_{k}^{-}< \tau_{k}^{(2)} <\tau_{k}^{+}  <\tau_{k}^{(3)}.$$

$$\tau_{k}^{-} \in \left[ \tau_{k}^{(1)}+ \dfrac{ \tau_{k}^{(2)}- \tau_{k}^{(1)}}{3} ,   \tau_{k}^{(2)}- \dfrac{ \tau_{k}^{(2)}- \tau_{k}^{(1)}}{3} \right]$$

and

$$\tau_{k}^{+} \in \left[ \tau_{k}^{(2)}+ \dfrac{ \tau_{k}^{(3)}- \tau_{k}^{(2)}}{3} ,   \tau_{k}^{(3)}- \dfrac{ \tau_{k}^{(3)}- \tau_{k}^{(2)}}{1} \right].$$

Also one has the equivalent
$$ |F_{k}(\tau_{k}^{\pm}) | \asymp   \dfrac{1}{ q_{n_{k}}(\alpha)^{2} q_{m_{k}}(\beta)^{2}}   |x_{k}-\tau_{k}^{\pm}|   |\tau_{k}^{\pm}-t_{k}(\alpha)||\tau_{k}^{\pm}-t_{k}(\beta)|    . $$

Given $\varepsilon >0$ arbitrary, replacing $F_{k}$ by $-F_{k}$  if necessary we are reduced to the three possible scenarios

\begin{enumerate}
\item  $\max\{ |F_{k}(\tau_{k}^{-}) |,  |F_{k}(\tau_{k}^{+}) | \}< \varepsilon $ then $J_{k}$ is a single interval.
\item $  F_{k}(\tau_{k}^{-}) < -\varepsilon$ and $  F_{k}(\tau_{k}^{+}) \leq \varepsilon$ $J_{k}$ is the union of two intervals.
\item $ -\varepsilon <  F_{k}(\tau_{k}^{-}) $ and $  \varepsilon< F_{k}(\tau_{k}^{+})$ $J_{k}$ is the union of two intervals.
\item $\min\{ |F_{k}(\tau_{k}^{-}) |,  |F_{k}(\tau_{k}^{+}) | \} \geq  \varepsilon$ $J_{k}$ is the union of three intervals.

\end{enumerate}

Unfortunately, we are not able to estimate the factors in $ |F_{k}(\tau_{k}^{\pm}) |$ due to the lack of a 2-dimensional analog of the notion of continued fraction for real numbers. We are going to treat the problem of the intervals by solving directly the equation $F_{k}(t)=\varepsilon$ in \S3.  We finish the study of the shape of the graph of $F_{k}$ by giving an upper bound on the measure of $t$ such that $|F_{k}(t)| \leq \varepsilon$.  This will follow from  the study of the set of the small values of the cubic polynomial  $F_{k}(t)$. A key result due to  H. Cartan \cite{c} given in the cubic case shows that the measure of set of small values cannot be too large. 
\begin{thm}[Theorem 3.1, VIII,\S 3 \cite{lg}] Let $f(t)=(t-t_{1})(t-t_{2})(t-t_{3})$ be a monic polynomial in $\mathbb{C}[t]$ and $\varepsilon$ be a positive real. Then,  the set of $t$ such that 
$$  |f(t)| \leq \varepsilon $$
is contained in the union of a most 3 intervals such that the sum of the lengths is bounded by $6e\varepsilon^{1/3}$. 
\end{thm}

\noindent Applying this result to the monic polynomial $p_{k}(t)=-A_{k}^{-1}F_{k}(t)$ tells us that there exists a covering of $\mathcal{D}(\varepsilon)$ consisting in at most three intervals $I_{1},I_{2}, I_{3} $ respectively with length $l_{1}, l_{2}, l_{3}$ so that their sums is less or equal to $6e|A_{k}| \varepsilon^{1/3}$. In particular, we have 
$$ \lambda (t \in \mathbb{R} ~:~ 0 < |F_{k}(t)| \leq \varepsilon) \leq \lambda(I_{1}^{k}) +  \lambda(I_{2}^{k}) +  \lambda(I_{3}^{k}) \leq  6e (|A_{k}^{-1}| \varepsilon)^{1/3}$$

$$  \lambda (t \in \mathbb{R} ~:~ 0 < |F_{k}(t)| \leq \varepsilon)  \ll (q_{n_{k}}(\alpha)^{2}  q_{m_{k}}(\beta)^{2} \varepsilon)^{1/3}. $$

\noindent We have obtained the following asymptotic estimate of the total length of the interval(s)

\begin{prop} \label{cartan}
For every positive real $\varepsilon$, there exists an integer $n_{0}$ such that  for all $k \geq k_{0}$,  
$$   \lambda (t \in \mathbb{R} ~:~ 0 < |F_{k}(t)| \leq \varepsilon)  \ll_{\varepsilon} (q_{n_{k}}(\alpha)  q_{m_{k}}(\beta))^{2/3}. $$
\end{prop}

\section{Intervals for the solutions of $F_{k}(t) = \pm \varepsilon$} \subsection*{Real solutions of the cubic equation $F_{k}(t) = \pm \varepsilon$}
Now we need to find all the possible intersection points which correspond to the values of $t$ such that $F_{k}(t)=\pm \varepsilon$. We have the expanded form 
$$F_{k}(t)= -A_{k} \left\lbrace t^{3}-  \Sigma_{1} t^{2} + \Sigma_{2}  t -  \Sigma_{3} \right\rbrace.$$

\noindent The equation $F_{k}(t)=\pm \varepsilon$ takes the following form

\begin{equation}
-A_k \{  t^{3} -\Sigma_{1}t^2 + \Sigma_{2}t-\Sigma_{3} \} \pm \varepsilon  = 0.
\end{equation}

\noindent We rather consider the equivalent equation

\begin{equation}\label{root}
G_{k}(\pm \varepsilon)(t) := t^{3} -\Sigma_{1}t^2 + \Sigma_{2}t-\left(\Sigma_{3} \pm \dfrac{\varepsilon}{A_k}\right) = 0.
\end{equation}
If we specializes $ \varepsilon $ to be zero, we recover the solution of the cubic equation $F_{k}(t)=0$. In particular we see that $F_{k}(t)$ and  $G_{k}(\pm \varepsilon)(t) $ share the same symmetric functions of the roots $ \Sigma_{1} $ and $\Sigma_{2} $, only the third differs by a factor  $ \pm \dfrac{\varepsilon}{A_k}$. \\
\noindent We choose the trigonometric parametrization of the solution of this cubic equation  $G_{k}(\pm \varepsilon)(t)=0$ following the presentation of \cite{bmc} (Appendix of Chapter 4).
The first step consists to perform the Vieta transform $t=y+\dfrac{\Sigma_{1}}{3}$ to get the following cubic in the reduced form (the $t^2$ term has vanished)
\begin{equation}\label{reduced}
y^{3}+P_{k}y=Q_{k}(\pm \varepsilon)
\end{equation}

\noindent where
\begin{center}
$P_{k}= \dfrac{3\Sigma_{2}-\Sigma_{1}^{2}}{3}$ and $Q_{k}(\pm \varepsilon)= \dfrac{-9\Sigma_{1}\Sigma_{2}+2\Sigma_{1}^{3} }{27}+\left( \Sigma_{3} \pm \dfrac{\varepsilon}{A_{k}} \right).$
\end{center}

\noindent The discriminant of $G_{k}(\pm \varepsilon)$ is given by 
$$  D_{k}(\pm \varepsilon) :=-4P_{k}^{3}-27Q_{k}(\pm \varepsilon)^{2}. $$

\noindent Let us solve (\ref{reduced}),  to do so we set $y= h_{k} z$ where $h_{k}= 2 \sqrt{|P_{k}|/3}$ and replace in (\ref{reduced})
$$h_{k}^{3} z^{3}+P_{k}h_{k} z  = \left( \dfrac{4|P_{k}|}{3} \right)^{3/2} z^3 + P_{k}  \left( \dfrac{4|P_{k}|}{3} \right) z = Q _{k}(\pm \varepsilon).$$

\noindent Then we divide by $ \left( \dfrac{4|P_{n}|}{3} \right)^{3/2}$

$$ z^3 +  \dfrac{3 P_{k}}{4|P_{k}|}  z =   \left( \dfrac{4|P_{k}|}{3} \right)^{-3/2} Q _{k}(\pm \varepsilon).$$

\noindent We are therefore reduced to solve the following equation, 

\begin{equation}\label{cubicreduced}
4 z^3 + 3~ \mathrm{sgn} (P_{k})  z = C_{k}(\pm \varepsilon)
\end{equation} 

\noindent where $$ C_{k}(\pm \varepsilon) =\dfrac{1}{2}  \left( \dfrac{3}{|P_{k}|} \right)^{3/2} Q _{k}(\pm \varepsilon). $$

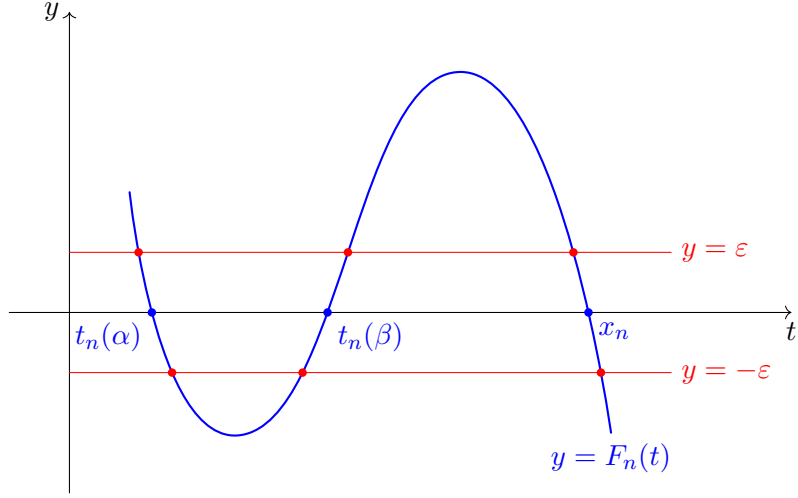
\begin{figure}

\begin{tikzpicture}[scale=0.8,cap=round]
 \tikzset{axes/.style={}}
 \begin{scope}[style=axes]
 \draw[->,name path=c4] (-1,0) -- (12,0) node[below] {$t$};
 \draw[->] (0,-3)-- (0,5) node[left] {$y$};
 \draw[red,name path=c1] (0,1)-- (10,1) node[right] {$y=\varepsilon$};
 \draw[red,-,name path=c2] (0,-1)-- (10,-1) node[right] {$y=-\varepsilon$};

 \draw [blue,thick,-,name path=c3] plot [smooth,tension=1] coordinates { 
   (1,2) (3,-2) (6.5,4) (9,-2)}  node[below] {$y=F_{k}(t)$}  ;
 
 \fill[red,name intersections={of=c1 and c3}]
    (intersection-1) circle (2pt)
    (intersection-2) circle (2pt)
        (intersection-3) circle (2pt) ;

        \fill[red,name intersections={of=c2 and c3}]
    (intersection-1) circle (2pt)
    (intersection-2) circle (2pt)
        (intersection-3) circle (2pt) ;

          \fill[blue,name intersections={of=c4 and c3}]
    (intersection-1) circle (2pt)node[below left]{$t_k(\alpha)$}
    (intersection-2) circle (2pt)node[below right]{$t_k(\beta)$}
        (intersection-3) circle (2pt) node[below right]{$x_k$};


\end{scope}
\end{tikzpicture}
\caption{The root $x_{k}$ dominates $ t_{k}(\alpha) $ and $t_{k}(\beta)$ as soon as $\delta_{k} >4/3$.}
\end{figure}

\noindent If we replace $ \varepsilon $ by zero, we recover the solution of the cubic equation $F_{k}(t)=0$. The latter equation has 3 real solutions $ x_{k} $, $ t_{k}(\alpha) $ and $ t_{k}(\beta) $, thus $ P_{k} <0 $ ( $P_k$ does not depend on $\varepsilon$). Therefore we are left with two  cases. \\
\noindent \textbf{Case 1.}$
\ \  \mathrm{sgn}(P_{k}) < 0, |C_{k}(\pm \varepsilon)| \geq 1$ in this case, one  real root one has to solve 
$$ 4 z^3 + 3  z = C_{k}(\pm \varepsilon).$$
Using the formula $\cosh 3 \theta = 4 \cosh^{3} \theta - 3\cosh^{2} \theta $

  $$z = \pm \mathrm{cosh} \left( \dfrac{1}{3} {\mathrm{arccosh}}(C_{k}(\pm \varepsilon))    \right)$$
  where the sign + is involved when $C_{k}(\pm \varepsilon) \geq 1$ and  the sign -  when $C_{k}(\pm \varepsilon) \leq -1$. If we assume that, say $C_{k}(\pm \varepsilon) \geq 1$  then we have a unique real root  for equation (\ref{reduced}) 
$$ y=2 \sqrt{-P_{k}/3} \cosh \left( \dfrac{1}{3}  \Phi_{k}(\pm \varepsilon) \right)$$
where $\Phi_{k}(\pm \varepsilon)=  \mathrm{arccosh} (C_{k}(\pm \varepsilon))$. Finally, we get the solutions of the equation $F_{k}(t)= \pm \varepsilon$
\begin{equation}\label{roots}
t_k(\pm \varepsilon) =2  \sqrt{-P_{k}/3}\cosh \left( \dfrac{1}{3}  \Phi_{k}(\pm \varepsilon) \right) + \dfrac{\Sigma_{1}}{3}.
\end{equation}

\noindent Thus $\mathcal{D}(\varepsilon)= \{ t \in \mathbb{R} , |F_{k}(t)|\leq \varepsilon \}$ is exactly the interval $I_{k} = [t_k({\varepsilon}), t_k({-\varepsilon})]$.

\noindent \textbf{Case 2.} $
 \ \  \mathrm{sgn}(P_{k}) < 0,  |C_{k}(\pm \varepsilon)| <1$. 

\noindent One is reduced to solve 
$$ 4 z^3 - 3  z = C_{k}(\pm \varepsilon).$$

\noindent Taking advantage of the relation $\cos 3 \theta = 4 \cos^{3} \theta - 3\cos^{2} \theta $  we deduce the solutions 

$$ z_{k}= \cos \left( \dfrac{1}{3}  \arccos (C_{k}(\pm \varepsilon)) + \dfrac{2(j-1)\pi}{3} \right)~~~~~(j=1,2,3).$$

\noindent Thus in the case when $|C_{k}(\pm \varepsilon)| \leq 1$, one has  three real distinct roots  for equation (\ref{reduced}) 
$$ y_{k}=2 \sqrt{-P_{n}/3} \cos \left( \dfrac{1}{3}  \Phi_{k}(\pm \varepsilon) + \dfrac{2(j-1)\pi}{3} \right)~~~~~(j=1,2,3)$$
where $\Phi_{k}(\pm \varepsilon)=  \arccos (C_{k}(\pm \varepsilon))$ lies in the open interval $(0,\pi)$. Finally, we get the solutions of the equation $F_{k}(t)= \pm \varepsilon$
\begin{equation}\label{roots}
t_{j,k}(\pm \varepsilon) =2  \sqrt{-P_{k}/3}\cos \left( \dfrac{1}{3}  \Phi_{k}(\pm \varepsilon) + \dfrac{2(j-1)\pi}{3} \right) + \dfrac{\Sigma_{1}}{3}~\  \ \ \ \ \ \ ~~~~(j=1,2,3).
\end{equation}

\noindent Thus $\mathcal{D}(\varepsilon)= \{ t \in \mathbb{R} , |F_{k}(t)|\leq \varepsilon \}$ is the union of three  the interval $I_{1}^{k} = [t_{1,k}({\varepsilon}), t_{1,k}({-\varepsilon})]$,  $I_{2}^{k} = [t_{2,k}({-\varepsilon}), t_{2,k}({\varepsilon})]$, $I_{3}^{k} = [t_{3,k}({\varepsilon}), t_{3,k}({-\varepsilon})]$. Hence we are in the worst case where $\mathcal{D}(\varepsilon)$ is the union of three intervals.

\subsection*{Asymptotic estimate for $ -P_{k} $} The next lemma is crucial it gives an asymptotic estimate for $-P_{k}$ with constants depending on $\varepsilon$.

\begin{lem}\label{estimatesPnQn}  Assume $\delta_{k} >4/3 $ for every $ k$, then as $ k $ gets large, one has\\

\begin{center}
$-P_{k} \asymp_{\varepsilon} q_{m_{k}}(\beta)^{2\delta_{k}}.$  
\end{center}

\end{lem}
\noindent \textit{Proof.} We have from Proposition \ref{sigmas}
 $$  q_{m_{k}}(\beta)^{\delta_{k}} + q_{m_{k}}(\beta)^{2-\delta_{k}} \ll_{\varepsilon} \Sigma_1  \ll_{\varepsilon} q_{m_{k}}(\beta)^{\delta_{k}} + q_{m_{k}}(\beta)^{2-\delta_{k}/2}.$$  

\noindent First, one has
$$ \delta_{k}-(2-\dfrac{\delta_{k}}{2})= \delta_{k}-2 +\dfrac{\delta_{k}}{2} = \dfrac{3\delta_{k}}{2}-2 >0.$$
 
\noindent  Also, 
 $$ \delta_{k}-(2-\delta_{k})= 2\delta_{k}-2  =2(\delta_{k}-1)  >0.$$
 
\noindent Thus, 
 $$  q_{m_{k}}(\beta)^{\delta_{k}}  \ll_{\varepsilon}  \Sigma_{1}  \ll_{\varepsilon} q_{m_{k}}(\beta)^{\delta_{k}}.$$
 
\noindent By definition $ -P_{k} = \dfrac{1}{3} \Sigma_{1}^{2}-\Sigma_{2}$ and it is positive because $ F_{k}(t) $ has three real roots. Therefore $ \Sigma_{1}^{2} $ is the dominant term of $ -P_{k} $, then 
$$  -P_{k} \asymp_{\varepsilon} \Sigma_{1}^{2} \asymp_{\varepsilon} q_{m_{k}}(\beta)^{2\delta_{k}}.$$

\section{Proof of Theorem \ref{main}}

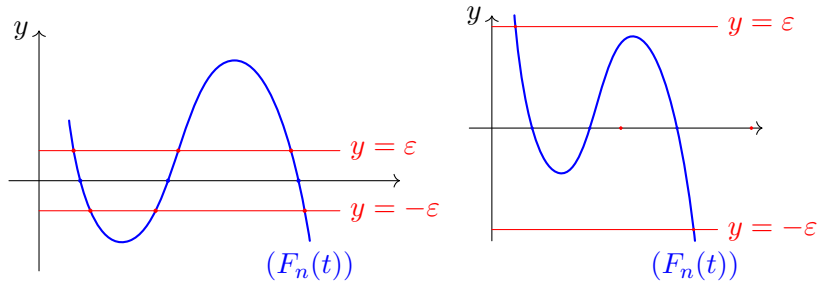
\begin{figure}
\begin{tikzpicture}[scale=0.4,cap=round]
 \tikzset{axes/.style={}}
 \begin{scope}[style=axes]
 \draw[->,name path=c4] (-1,0) -- (12,0) node[below] {};
 \draw[->] (0,-3)-- (0,5) node[left] {$y$};
 \draw[red,name path=c1] (0,1)-- (10,1) node[right] {$y=\varepsilon$};
 \draw[red,-,name path=c2] (0,-1)-- (10,-1) node[right] {$y=-\varepsilon$};

 \draw [blue,thick,-,name path=c3] plot [smooth,tension=1] coordinates { 
   (1,2) (3,-2) (6.5,4) (9,-2)}  node[below] {$(F_{n}(t))$}  ;
 
 \fill[red,name intersections={of=c1 and c3}]
    (intersection-1) circle (2pt)
    (intersection-2) circle (2pt)
        (intersection-3) circle (2pt) ;

        \fill[red,name intersections={of=c2 and c3}]
    (intersection-1) circle (2pt)
    (intersection-2) circle (2pt)
        (intersection-3) circle (2pt) ;

          \fill[blue,name intersections={of=c4 and c3}]
    (intersection-1) circle (2pt)node[below left]{}
    (intersection-2) circle (2pt)node[below right]{}
        (intersection-3) circle (2pt) node[below right]{};


\end{scope}
\end{tikzpicture}
\begin{tikzpicture}[scale=0.3,cap=round]
 \tikzset{axes/.style={}}
 \begin{scope}[style=axes]
 \draw[->,name path=c4] (-1,0) -- (12,0) node[below] {};
 \draw[->] (0,-5)-- (0,5) node[left] {$y$};
 \draw[red,name path=c1] (0,4.5)-- (10,4.5) node[right] {$y=\varepsilon$};
 \draw[red,-,name path=c2] (0,-4.5)-- (10,-4.5) node[right] {$y=-\varepsilon$};

 \draw [blue,thick,-,name path=c3] plot [smooth,tension=1] coordinates { 
   (1,5) (3,-2) (6.5,4) (9,-5)}  node[below] {$(F_{n}(t))$}  ;
 
 \fill[red,name intersections={of=c1 and c3}]
    (intersection-1) circle (2pt)
    (intersection-2) circle (2pt)
        (intersection-3) circle (2pt) ;

        \fill[red,name intersections={of=c2 and c3}]
    (intersection-1) circle (2pt)
    (intersection-2) circle (2pt)
        (intersection-3) circle (2pt) ;

          \fill[blue,name intersections={of=c4 and c3}]
    (intersection-1) circle (2pt)node[below left]{}
    (intersection-2) circle (2pt)node[below right]{}
        (intersection-3) circle (2pt) node[below right]{};


\end{scope}
\end{tikzpicture}
\caption{One the left the worst scenario $|C_{k}(\pm \varepsilon)| <1$  (three intervals) and one the right the best scenario $|C_{k}(\pm \varepsilon)| \geq 1$ (one interval).}
\end{figure}

\noindent Let us fix $\varepsilon >0$, and $ N_{k}=q_{m_{k}}(\beta)^{\delta_{k}} $ where $\delta_{k} > 4/3$ is chosen such that 
\begin{equation}\label{choicedelta}
 \delta_{k} \in \left( \dfrac{4}{3}, \dfrac{1}{5}(3+4\gamma_{k}-\eta_{k}) \right]\neq \emptyset
\end{equation}
for $k$ large enough. The nonemptyness in (\ref{choicedelta}) is ensured by condition (\ref{gammabound}), indeed the length of  $\left( \dfrac{4}{3}, \dfrac{1}{5}(3+4\gamma_{k}-\eta_{k}) \right]$  is given by
$$\dfrac{1}{5}(3+4\gamma_{k}-\eta_{k})-\dfrac{4}{3} =\dfrac{4}{5}\gamma_{n}-\dfrac{1}{5}\eta_{k}-\dfrac{11}{15}> \dfrac{4}{5}(\dfrac{11}{12}+\dfrac{\eta_{k}}{4})-\dfrac{1}{5}\eta_{k}-\dfrac{11}{15}=0.$$
\noindent We assume we are in the worst case\footnote{The shape of $ \mathcal{D}(\varepsilon) $ can convince the reader that the best case is when the line cuts  he boundary of $ \mathcal{D}(\varepsilon) $ only twice.}, namely the case 2 when $ |C_{k}(\pm \varepsilon)|<1$. Thus $F_{k}(t) = \pm \varepsilon$ has three real roots. In that case we have seen in \ref{roots} that the three (trigonometric) solutions are given by 
$$ t_k(\pm \varepsilon) =2  \sqrt{-P_{k}/3}\cos \left( \dfrac{1}{3}  \Phi_{k}(\pm \varepsilon) + \dfrac{2(j-1)\pi}{3} \right) + \dfrac{\Sigma_{1}}{3}~\  \ \ \ \ \ \ ~~~~(j=1,2,3). $$
\noindent Thus, the set of times $t$ such that $|F_{k}(t)| \leqslant \varepsilon$ is the union of three intervals $I_{1}^{k}$, $I_{2}^{k}$ and $I_{3}^{k}$ bounded respectively by the roots $t_{1}(\pm \varepsilon),t_{2}(\pm \varepsilon)$ and $t_{3}(\pm \varepsilon)$. 

 \noindent Since the lengths of each interval $ I_{1}^{k} $,$ I_{2}^{k} $ and $ I_{3}^{k} $ are all equal, we only need to focus on the length of $I_{1}^{k}$.  One has
 $$ t_1(\pm \varepsilon) =2  \sqrt{-P_{k}/3}\cos \left( \dfrac{1}{3}  \Phi_{k}(\pm \varepsilon) \right) + \dfrac{\Sigma_{1}}{3}. $$

\noindent The length of $I_{1}^{k}$ is given by 
$$ \lambda (I_{1}^{k})= \mathrm{length}(I_{1}^{k}) = |t_{1}(\varepsilon) -t_{1}(-\varepsilon)| =2  \sqrt{-P_{k}/3}\left| \cos \dfrac{\Phi_{k}(\varepsilon)}{3} - \cos \dfrac{\Phi_{k}(-\varepsilon)}{3} \right|.$$
Bounds for the length of $I_{1}^{k}$  are given by, 
\begin{prop}\label{lengthinterval} For $n$ large enough, 
$$     q_{m_{k}}(\beta)^{4(1+ \gamma_{k})-5\delta_{k}} \ll_{\varepsilon} \lambda (I_{1}^{k}) \ll_{\varepsilon}   q_{m_{k}}(\beta)^{2/3(1+\gamma_{k})}.$$
\end{prop}

\noindent \textit{Proof of the Proposition.} The upper bound is a direct consequence of Proposition \ref{cartan}. Now let us estimate the lower bound.
\noindent Using the identity $\displaystyle \cos \alpha - \cos \beta = -2 \sin(\frac{\alpha+ \beta}{2})\sin(\frac{\alpha- \beta}{2})$,  we get
$$  \lambda (I_{1}^{k})=4  \sqrt{-P_{k}/3} \left| \sin\left( \dfrac{\Phi_{k}(\varepsilon) + \Phi_{k}(-\varepsilon)}{6}\right)  \sin \left( \dfrac{\Phi_{k}(\varepsilon) - \Phi_{k}(-\varepsilon)}{6} \right) \right|. $$

\noindent Remark that, one has $|(\Phi_{k}(-\varepsilon) \pm \Phi_{k}(\varepsilon))/6 |\in [0,\pi/3]$, and thus one can apply the inequality $\sin x \geqslant x/2$ which is valid for all $x \in [0,\pi/2]$.  We obtain
$$ \lambda (I_{1}^{k}) \geq \dfrac{1}{36} \sqrt{-P_{k}/3}~ 
\left|\Phi_{k}(\varepsilon)^{2} -\Phi_{k}(-\varepsilon)^{2} \right|.  $$
or more precisely, 
$$ \lambda (I_{1}^{k}) \geq \dfrac{1}{36} \sqrt{-P_{k}/3}~ 
\left|\arccos^{2}(C_{k}(\varepsilon))  -\arccos^{2}(C_{k}(-\varepsilon)) \right|.  $$

\begin{lem} For all $ x,y \in [-1,1] $, 
$$ | \arccos^{2} x -  \arccos^{2} y | \geq |x-y|^{2}. $$

\end{lem}
\noindent \textit{Proof.} One remarks that 
$$ \left| \arccos^{2} x - \arccos^{2} y  \right| =  \left| \arccos x - \arccos y  \right| | \arccos x + \arccos y |.  $$
\noindent Noting the fact that the function arcos assumes positive values on the open interval $(-1,1)$, we can write that
$$ | \arccos x + \arccos y |=\arccos x + \arccos y  \geq  | \arccos x - \arccos y|.$$
\noindent Thus, 
$$ \left| \arccos^{2} x - \arccos^{2} y  \right| \geq  \left| \arccos x - \arccos y  \right|^{2}.  $$
\noindent Moreover, the Lipschitz property for the cosine yields
$$|\cos u -\cos v|  \leq  | u -v |.$$
The latter applied to $u=\arccos x$ and $v=\arccos y$ we get 
\begin{equation}\label{ineqarcos}
|x-y|  \leq  |\arccos x - \arccos x |.
\end{equation}

\noindent Hence, 
$$ \left| \arccos^{2} x - \arccos^{2} y  \right| \geq  \left| x - y  \right|^{2}.  $$

\begin{flushright}
$  \square$
\end{flushright}

\noindent Remembering that  $\Phi_{k}(\pm \varepsilon)=  \arccos (C_{k}(\pm \varepsilon))$ the previous lemma gives us
$$|C_{k}( \varepsilon)-C_{k}(- \varepsilon)|^{2}  \leq  |\Phi_{k}(\varepsilon)^{2}  -\Phi_{k}(- \varepsilon)^{2} |.$$

\noindent Therefore, 
$$  \lambda (I_{1}^{k}) \geq \dfrac{1}{36 \sqrt{3}} \sqrt{-P_{k}}~\left|\Phi_{k}(\varepsilon)^{2} -\Phi_{k}(-\varepsilon)^{2} \right| \gg  \sqrt{-P_{k}} \left|  C_{k}(\varepsilon)-C_{k}(-\varepsilon)\right|^{2}.  $$

\noindent Thus, we obtain
$$ \lambda (I_{1}^{k}) \gg \sqrt{-P_{k}}\left|  C_{k}(\varepsilon)-C_{k}(-\varepsilon) \right|^{2}.  $$

\noindent Concerning the right hand side, one has
$$ C_{k}(\varepsilon)-C_{k}(-\varepsilon) = \dfrac{1}{2}  \left( \dfrac{3}{-P_{k}} \right)^{3/2} \left(  Q _{k}( \varepsilon) -   Q _{k}(-\varepsilon)   \right)=  \dfrac{1}{2}  \left( \dfrac{3}{-P_{k}} \right)^{3/2}\dfrac{2\varepsilon}{A_{k}}  = \dfrac{3\sqrt{3} \varepsilon}{(-P_{k})^{3/2}A_{k}}.$$

\noindent Thus, for $n$ large enough
$$ |C_{k}(\varepsilon)-C_{k}(-\varepsilon)| \asymp \dfrac{  q_{m_{k}}(\beta)^{2(1+ \gamma_{k})} }{(-P_{k})^{3/2}} $$

\noindent since $ A_{k}^{-1}= (e_{n_{k}}(\alpha)e_{m_{k}}(\beta))^{-1} \sim q_{n_{k}}(\alpha)^{2}q_{m_{k}}(\beta)^{2} = q_{m_{k}}(\beta)^{2(1+\gamma_{k})}$ by (\ref{error2m}). Hence, we get the lower bound for the length

 \begin{equation}
  \lambda (I_{1}^{k}) \gg_\varepsilon  \sqrt{-P_{k}}\left|  C_{k}(\varepsilon)-C_{k}(-\varepsilon) \right|^{2}  \gg_\varepsilon   \dfrac{  q_{m_{k}}(\beta)^{4(1+ \gamma_{k})} }{(-P_{k})^{5/2}} .
 \end{equation}

\begin{flushright}
$ \square $
\end{flushright}

\noindent The lemma \ref{estimatesPnQn} yields the expected lower bound,

 \begin{equation}
  \lambda (I_{1}^{k}) \gg_\varepsilon    q_{m_{k}}(\beta)^{4(1+ \gamma_{k}) -5\delta_{k}} .
 \end{equation}

\begin{flushright}
$  \square$
\end{flushright}

\noindent We are ready to complete the proof of the theorem. One hypothesis of the theorem tells us that the least common multiple $ l_{k} $ of $q_{n_{k}}(\alpha)  $  and $q_{m_{k}}(\beta)  $ satisfies 
$$  l_{k} \leq q_{m_{k}}(\beta)^{1+\eta_{k}}$$
with some $0\leq \eta_{k} \leq 1/3  $.
Using the latter with Lemma \ref{lengthinterval}, 

\begin{equation}
  \dfrac{ \lambda (I_{1}^{k})}{l_{k}} \gg_\varepsilon  \dfrac{  q_{m_{k}}(\beta)^{4(1+ \gamma_{k})-5\delta_{k}} }{q_{m_{k}}(\beta)^{1+ \eta_{k}}}
 \end{equation}

\noindent which reduces to 
\begin{equation}
  \dfrac{ \lambda (I_{1}^{k})}{l_{k}} \gg_\varepsilon    q_{m_{k}}(\beta)^{3+4\gamma_{k}-\eta_{k}-5\delta_{k}}.
 \end{equation}

  \noindent For each $k$, let us consider $ j_{k} \in \{1,2,3\} $ to be the index of  one of two intervals $ I_{j}^{k} $ which do not contain $ x_{k} $. The initial choice of the sequence $\delta_{k}$ in (\ref{choicedelta}) implies that $ 3+4\gamma_{k}-\eta_{k}-5\delta_{k} $ is positive for $n$ large enough. In particular, 
  \begin{equation}\label{divergencelength}
 \lim_{k} \dfrac{ \lambda (I_{1}^{k})}{l_{k}} =\infty.
 \end{equation}

\noindent   The same holds for  $\lambda( I_{j_{k}}^{k})=  \lambda (I_{1}^{k})$ and therefore $I_{j_{k}}^{k}$ always contains at least one multiple of $ l_{k}$. For each $k$ large enough, let $ t_{k} = a_{k} l_{k} \in I_{j_{k}}$ be a multiple of $l_k$, and by definition of $j_{k}$,  $ x_{k} \neq t_{k} $. Moreover, if we set $ u_{k} =v_{k}(t_{k}) $, we have
  $$u_{k}=(x_{k}-t_{k}, y_{k}-a_{k} \lambda_{k} c_{n_{k}}(\alpha), z_{k}-a_{k} \lambda_{k}c_{m_{k}}(\beta)  ) \in \mathbb{Z}^{3} \cap \mathcal{D}(\varepsilon).  $$ 

\noindent In other words, we have found $u_{k} \in \mathbb{Z}^{3}$ with nonzero first component which satisfies 
$$ 0 < |f(u_{k})| \leq \varepsilon. $$
In particular, 
$$ \inf_{u\in \mathbb{Z}^{3}, u_{1}\neq 0} |f(u)|=0.$$
Hence $(\alpha,\beta)$ satisfies the Littlewood conjecture. The theorem \ref{main} is proved.
\begin{flushright}
$  \square$
\end{flushright}


\begin{thebibliography}{ZZZZ}

\bibitem{ab} Adamczewski  B. and Bugeaud Y.
    On the Littlewood conjecture in simultaneous Diophantine approximation,  J. London Math. Soc. 73 (2006) 355–366.

  \bibitem{bmc} Birkhoff G. and S. Mac Lane, Survey of Modern Algebra, CRC Press 4th Ed.   1 (1997).
  
  
\bibitem{b}  Bugeaud Y.
    Around the Littlewood conjecture in Diophantine approximation,  Publ. Math. Bes. 1 (2014) 5-18.
    
 \bibitem{c}  Cartan H. \textit{Sur les syst\`emes de fonctions holomorphes \`a vari\'et\'es lin\'eaires lacunaires et leurs applications}, Annales scientifiques de l'\'Ecole Normale Sup\'erieure, 3e s{\'e}rie, 45, 255-346 (1928).    

\bibitem{cs} Cassels J. W. S. and Swinnerton-Dyer H. P. F. 
    On the product of three homogeneous linear forms and the indefinite ternary quadratic forms  Philos. Trans. Roy. Soc. London. Ser. A. 248 (1955) 73–96.
    
 \bibitem{dm}   De Mathan B. 
   \textit{ Conjecture de Littlewood et récurrences linéaires }
    J. Théor. Nombres Bordeaux 15 (2003) 249–266.
    
    
     
\bibitem{ekl}    Einsiedler M. and Katok A. and Lindenstrauss E. Invariant measures and the set of exceptions to Littlewood's conjecture  Ann. of Math. 164 (2006) 513–560.

\bibitem{ew}    Einsiedler M. and Thomas W. Ergodic Theory with a view towards Number Theory, GTM  259, Springer 2011. 

\bibitem{hw}   Hardy G.H.  and Wright E.M.  An introduction to the theory of numbers, 6th Ed. (extended version by Heath-Brown R. and Silverman J.)  Oxford University Press 2008. 



\bibitem{little}  Littlewood  J.E. , Some Problems in Real and Complex analysis, D.C. Heath; First Ed. (1968).
\bibitem{lg}    Lang S. Introduction to Complex Hyperbolic Spaces, Springer 1987. 

\bibitem{l} Lazar Y. \textit{Explicit solutions to the Oppenheim conjecture for indefinite ternary diagonal forms}~ Acta Arith., 205.4, (2022), 287-307.


\bibitem{m}  Mahler K.  \textit{On the continued fractions of quadratic and cubic irrationals}, Ann. Mat. Pura Appl. (4) 30 (1949), 147-172. 
    
 \bibitem{m1} Margulis G.A.  \textit{Indefinite quadratic forms and unipotent flows on homogeneous spaces,} Banach Center Publish. vol. 23, Polish Scientific Publishers, Warsaw, 1989.


\bibitem{m2}~Margulis G.A.  \textit{Oppenheim Conjecture} in Fields Medalist Lectures, World Scientific Series in the 20th Century,  9, Sec. Ed.  (2003), 281-336.
  
   \bibitem{pv}      Pollington A. D. and Velani S. 
    On a problem in simultaneous Diophantine approximation: Littlewood's conjecture  Acta Math. 185 (2000) 287–306.
  

   \bibitem{q1}  Queff\'{e}lec M. \textit{Transcendance des fractions continues de Thue-Morse} , J. Number Theory   73 (1998), 201 - 211.
   
   \bibitem{q2}  Queff\'{e}lec M. Une introduction \`{a} la conjecture de Littlewood, S\'{e}minaires et Congr\`{e}s 20, (2009), p. 129 - 152.
   
   \bibitem{rs} Rahman Q.R. ~and ~Schmeisser G. Analytic Theory of polynomials,  London Math. Soc. Monographs, New series \textbf{26}, Oxford University Press, Press Inc. New-York  2002. 
   
   
   
   
\bibitem{s}
Schmidt Wolfgang~M. 
\newblock {\em Diophantine approximation}, volume 785 of {\em Lecture Notes in
  Mathematics}.
\newblock Springer, Berlin, 1980.


 \bibitem{v}   Venkatesh  A.
    The work of Einsiedler, Katok and Lindenstrauss on the Littlewood conjecture 
    Bull. Amer. Math. Soc. (N.S.) 45 (2008) 117–134.





\end{thebibliography}
\end{document}